\documentclass[10pt]{amsart}
\usepackage{amsfonts}
\usepackage{graphicx}
\usepackage{amscd}
\usepackage{amsmath,amsfonts,amssymb,amsthm,mathrsfs,xypic}
\xyoption{curve}
\usepackage{fancybox}
\newtheorem{thm}{Theorem}[section]

\newtheorem{cor}[thm]{Corollary}

\newtheorem{lem}[thm]{Lemma}

\newtheorem{prop}[thm]{Proposition}
\newtheorem{rem}[thm]{Remark}
\newtheorem{defn}[thm]{Definition}
\numberwithin{equation}{section}

\newcommand{\m}{\Lambda}

\newcommand{\s}{\hfill\blacksquare}

\newcommand{\Cok}{\operatorname{Coker}}
\newcommand{\Ker}{\operatorname{Ker}}
\newcommand{\Hom}{\operatorname{Hom}}

\begin{document}
\title [Symmetric recollements]{Symmetric recollements induced by bimodule extensions}
\author [Pu Zhang] {Pu Zhang}
\thanks{Supported by the NSF of China (10725104), and STCSM (09XD1402500).}
\thanks{pzhang$\symbol{64}$sjtu.edu.cn}

\maketitle

\centerline{Department of Mathematics, \ \ Shanghai Jiao Tong
University} \centerline{Shanghai 200240, P. R. China}
\begin{abstract} Inspired by the work of J$\o$rgensen [J], we define a
(upper-, lower-) symmetric recollements; and give a one-one
correspondence between the equivalent classes of the upper-symmetric
recollements and one of the lower-symmetric recollements, of a
triangulated category. Let $\m = \left(\begin{smallmatrix}
A&M\\
0&B
\end{smallmatrix}\right)$ with  bimodule $_AM_B$.
We construct an upper-symmetric abelian category recollement of
$\m$-mod; and a symmetric triangulated category recollement of
$\underline {\m\mbox{-}\mathcal Gproj}$ if $A$ and $B$ are
Gorenstein and $_AM$ and $M_B$ are projective.

\vskip5pt

\end{abstract}

\vskip10pt

\centerline {\bf  Introduction}

\vskip10pt

A triangulated category recollement, introduced by A. A. Beilinson,
J. Bernstein, and P. Deligne [BBD], and an abelian category
recollement, formulated by V. Franjou and T. Pirashvili [FV], play
an important role in algebraic geometry and in representation theory
([MV], [CPS], [K], [M]).

\vskip5pt

Recently, P. J$\o$rgensen [J] observed that if a triangulated
category $\mathcal C$ has a Serre functor, then a triangulated
category recollement of
 $\mathcal C$ relative to $\mathcal C'$ and
$\mathcal C''$ can be interchanged in two ways to triangulated
category recollements of
 $\mathcal C$ relative to $\mathcal C''$ and
$\mathcal C'$. Inspired by [J] we define in Section 2 a (upper-,
lower-) symmetric recollement; and prove that there is a one-one
correspondence between the equivalent classes of the upper-symmetric
triangulated category recollements of $\mathcal C$ relative to
$\mathcal C'$ and $\mathcal C''$, and the ones of the
lower-symmetric triangulated category recollements of $\mathcal C$
relative to $\mathcal C''$ and $\mathcal C'$. Let $A$ and $B$ be
Artin algebras, $M$ an $A$-$B$-bimodule, and $\m =
\left(\begin{smallmatrix}
A&M\\
0&B
\end{smallmatrix}\right)$ the upper triangular matrix algebra.
We construct an upper-symmetric abelian category recollement of
$\m$-mod, the category of finitely generated $\m$-modules. \vskip5pt

An important feature of Gorenstein-projective modules is that the
category $A\mbox{-}\mathcal Gproj$ of Gorenstein-projective
$A$-modules is a Frobenius category, and hence the stable category
$\underline {A\mbox{-}\mathcal Gproj}$ is a triangulated category
([Hap]). Iyama-Kato-Miyachi ([IKM], Theorem 3.8) prove that if $A$
is a Gorenstein algebra, then $\underline {T_2(A)\mbox{-}\mathcal
Gproj}$ admits a triangulated category recollement, where $T_2(A) =
\left(\begin{smallmatrix}
A&A\\
0&A
\end{smallmatrix}\right)$. In Section 3, if $A$ and $B$
are Gorenstein algebras and $_AM$ and $M_B$ are projective, we
extend this result by asserting that \ $\underline
{\m\mbox{-}\mathcal Gproj}$ admits a symmetric triangulated category
recollement, and by explicitly writing out the involved functors.

\section{\bf \ An equivalent definition of triangulated category recollements}

\subsection{} \ Recall the following

\begin{defn} \ $(1)$ \ ([BBD]) \ Let ${\mathcal C}', \ {\mathcal C}$  and
${\mathcal C}^{\prime\prime}$ be  triangulated categories. The
diagram

\begin{center}
\begin{picture}(120,45)
\put(0,20){\makebox(3,1){${\mathcal C}'$}}
\put(40,30){\vector(-1,0){30}} \put(10,20){\vector(1,0){30}}
\put(40,10){\vector(-1,0){30}} \put(50,20){\makebox(3,1){${\mathcal
C}$}} \put(90,30){\vector(-1,0){30}} \put(60,20){\vector(1,0){30}}
\put(90,10){\vector(-1,0){30}} \put(100,20){\makebox(3,1){${\mathcal
C}^{\prime\prime}$}} \put(25,35){\makebox(3,1){\scriptsize$i^\ast$}}
\put(25,25){\makebox(3,1){\scriptsize$i_\ast$}}
\put(25,5){\makebox(3,1){\scriptsize$i^!$}}
\put(75,35){\makebox(3,1){\scriptsize$j_!$}}
\put(75,25){\makebox(3,1){\scriptsize$j^\ast$}}
\put(75,5){\makebox(3,1){\scriptsize$j_\ast$}}
\put(240,15){\makebox(25,1){$(1.1)$}}
\end{picture}
\end{center}
\noindent of {\bf exact} functors is a triangulated category
recollement of  $\mathcal C$ relative to $\mathcal C'$ and $\mathcal
C''$, if the following conditions are satisfied:

\vskip5pt

${\rm (R1)}$ \  $(i^\ast, i_\ast)$,  \ $(i_\ast, i^!)$, \ $(j_!,
j^\ast)$, and $(j^\ast, j_\ast)$ are adjoint pairs;

\vskip5pt

${\rm (R2)}$  \ $i_*,  \ j_!$ and $j_*$ are fully faithful;

\vskip5pt

${\rm (R3)}$ \ $j^*i_\ast = 0$;

\vskip5pt

${\rm (R4)}$ \ For each object $X\in \mathcal C$, the counits and
units  give rise to distinguished triangles:
$$j_!j^*X \stackrel {\epsilon_{_X}}\longrightarrow X \stackrel {\eta_{_X}}\longrightarrow i_*i^*X\longrightarrow \  \ \ \
\mbox{and} \ \ \ \ i_*i^!X \stackrel {\omega_{_X}}\longrightarrow
X\stackrel {\zeta_{_X}}\longrightarrow j_*j^*X\longrightarrow .$$

\vskip10pt

$(2)$ \ ([FV]) \ Let ${\mathcal C}', \ {\mathcal C}$  and ${\mathcal
C}^{\prime\prime}$ be  abelian categories. The diagram $(1.1)$ of
{\bf additive} functors is an abelian category recollement of
$\mathcal C$ relative to $\mathcal C'$ and $\mathcal C''$,  if \
${\rm (R1)}$, ${\rm (R2)}$ and ${\rm (R5)}$ are satisfied, where

\vskip5pt

 ${\rm (R5)}$ \ ${\rm Im}i_\ast = {\rm Ker}j^*$.
\end{defn}

\vskip10pt

\begin {rem} \ $(1)$ \ Let $(1.1)$ be an abelian category recollement. If all the
involved functors are exact, then one can prove that there is an
equivalence $\mathcal C\cong \mathcal C'\times \mathcal C''$ of
categories. This explains why Franjou-Pirashvili [FV] did not
require the exactness of the involved functors in Definition
1.1$(2)$.

\vskip10pt

$(2)$ \ For any adjoint pair $(F, G)$, it is well-known that $F$ is
fully faithful if and only if the unit $\eta: {\rm Id} \rightarrow
GF$ is a natural isomorphism, and $G$ is fully faithful if and only
if the counit $\epsilon: FG\longrightarrow {\rm Id}$ is a natural
isomorphism; and that if $F$ is fully faithful then $G\epsilon_{_X}$
is an isomorphism for each object $X$, and if $G$ is fully faithful
then $F\eta_{_X}$ is an isomorphism for each object $Y$.

\vskip10pt

$(3)$ \  In any triangulated or abelian category recollement, under
the condition ${\rm (R1)}$, the condition ${\rm (R2)}$ is equivalent
to the condition ${\rm (R2')}$: the units ${\rm Id}_{\mathcal
C'}\rightarrow i^!i_\ast$ and ${\rm Id}_{{\mathcal
C}^{\prime\prime}}\rightarrow j^\ast j_!$, and the counits $i^\ast
i_\ast\rightarrow {\rm Id}_{{\mathcal C}^{\prime}}$ and $j^\ast
j_\ast\rightarrow {\rm Id}_{{\mathcal C}^{\prime\prime}}$, are
natural isomorphisms.

\vskip10pt

$(4)$ \ In an abelian category recollement one has $i^*j_! = 0$ and
$i^!j_* = 0$; and in a triangulated category recollement one has
${\rm Im}i_\ast = {\rm Ker}j^*$, ${\rm Im}j^! = {\rm Ker}i^*$ and \
${\rm Im}j_\ast = {\rm Ker}i^!$.

\vskip10pt

$(5)$ \ In any abelian category recollement $(1.1)$, the counits and
units give rise to exact sequences of natural transformations \
$j_!j^* \rightarrow {\rm Id}_{\mathcal A}\rightarrow
i_*i^*\rightarrow 0$ and $0\rightarrow i_*i^! \rightarrow {\rm
Id}_{\mathcal A}\rightarrow j_*j^*;$ and if $\mathcal C', \mathcal
C$, and $\mathcal C''$ have enough projective objects, then $i^*$ is
exact if and only if $i^!j_! = 0$; and dually, if $\mathcal C',
\mathcal C$, and $\mathcal C''$ have enough injective objects, then
$i^!$ is exact if and only if $i^*j_* = 0$. See [FV].
\end{rem}

\subsection{} \ We will need the following equivalent definition of a triangulated category recollement, which
possibly makes the
construction of a triangulated category recollement easier.

\begin{lem} Let $(1.1)$ be a diagram of exact functors of triangulated categories. Then
it is a triangulated category recollement if and only if the
conditions ${\rm (R1)}$, ${\rm (R2)}$ and ${\rm (R5)}$ are
satisfied.
\end{lem}

\noindent{\bf Proof.} \ This seems to be well-known, however we did
not find an exact reference. For the convenience of the reader we
include a proof.

We only need to prove the sufficiency. Embedding the counit morphism
$\epsilon_{_X}$ into a distinguished triangle $j_!j^*X \stackrel
{\epsilon_{_X}}\rightarrow X\stackrel h\rightarrow Z\rightarrow$.
Applying $j^*$ we get a distinguished triangle $j^*j_!j^*X \stackrel
{j^*\epsilon_{_X}} \rightarrow j^*X \stackrel {j^*h}\rightarrow
j^*Z\rightarrow$. Since $j^*\epsilon_{_X}$ is an isomorphism by
Remark 1.2$(2)$, we have $j^*Z = 0$. By ${\rm Im}i_\ast = {\rm
Ker}j^*$ we have $Z = i_*Z'$. Applying $i^*$ to the distinguished
triangle $j_!j^*X \stackrel {\epsilon_{_X}}\rightarrow X\stackrel
h\rightarrow i_*Z'\rightarrow$,  by $i^*j_! = 0$ we know that $i^*h:
\ i^*X \rightarrow i^*i_*Z'$ is an isomorphism. Since the counit
morphism $i^*i_*Z' \stackrel {\varepsilon_{_{Z'}}} \longrightarrow
Z'$ is an isomorphism, we have isomorphism
$i_*((i^*h)^{-1})i_*(\varepsilon_{_{Z'}}^{-1}): \ i_*Z'\rightarrow
i_*i^*X$, and hence we get a distinguished triangle of the form
$j_!j^*X \stackrel {\epsilon_{_X}}\rightarrow X\stackrel
f\rightarrow i_*i^*X\rightarrow$ \   with $f =
i_*((i^*h)^{-1})i_*(\varepsilon_{_{Z'}}^{-1})h$, which also means
${\rm Im}j_! = {\rm Ker}i^*$. Since $i^*h$ is an isomorphism, $i^*f$
is an isomorphism.  \vskip5pt

In order to complete the first distinguished triangle in ${\rm
(R4)}$, we need to show that $f$ can be chosen to be the unit
morphism. Embedding the unit morphism $\eta_{_X}$ into a
distinguished triangle $Y \rightarrow X \stackrel
{\eta_{_X}}\rightarrow i_*i^*X\rightarrow$. By the similar argument
(but this time we use ${\rm Im}j_! = {\rm Ker}i^*$) we get a
distinguished triangle of the form $j_!j^*X \stackrel g\rightarrow X
\stackrel {\eta_{_X}}\rightarrow i_*i^*X\rightarrow.$ By the
following commutative diagram given by the adjoint pair $(i^*, i_*)$
\[
\xymatrix{\Hom_{\mathcal C'}(i^*i_*i^*X, i^*X)
\ar[r]^-{\sim}\ar[d]^{(i^*f, -)}_-{\wr}  & \Hom_{\mathcal
C}(i_*i^*X, i_*i^*X)
\ar[d]^{(f, -)}\\
 \Hom_{\mathcal C'}(i^*X, i^*X) \ar[r]^-{\sim}  & \Hom_{\mathcal C}(X, i_*i^*X)}
\]
we see that $\Hom_{\mathcal C}(f, i_*i^*X)$ is also an isomorphism,
and hence there is $u\in\Hom_{\mathcal C}(i_*i^*X, i_*i^*X)$ such
that $uf = \eta_{_X}$. Since $(i^*, i_*)$ is an adjoint pair and
$i_*$ is fully faithful, it follows that $i^*\eta_{_X}$ is an
isomorphism. Replacing $f$ by $\eta_{_X}$ we get $v\in\Hom_{\mathcal
C}(i_*i^*X, i_*i^*X)$ such that $v\eta_{_X} = f$. Thus we have
morphisms of distinguished triangles
\[\xymatrix {
j_!j^*X \ar[r]^{\epsilon_{_X}}\ar[d]^-{=} & X\ar[d]^-{=}\ar[r]^-{f} & i_*i^*X\ar[d]_{vu}\ar[r]& \\
j_!j^*X \ar[r]^{\epsilon_{_X}}& X\ar[r]^-{f} & i_*i^*X\ar[r]&}\] and
\[\xymatrix {
j_!j^*X \ar[r]^{g}\ar[d]^-{=} & X\ar[d]^-{=}\ar[r]^-{\eta_{_X}} & i_*i^*X\ar[d]_{uv}\ar[r]&\\
j_!j^*X \ar[r]^{g}& X\ar[r]^-{\eta_{_X}} & i_*i^*X\ar[r]&.}\] So
$uv$ and $vu$, and hence $u$ and $v$, are isomorphisms. By the
isomorphism of triangles
\[\xymatrix {
j_!j^*X \ar[r]^{\epsilon_{_X}}\ar[d]^-{=} & X\ar[d]^-{=}\ar[r]^-{\eta_{_X}} & i_*i^*X\ar[d]^{\wr}_{v}\ar[r]& \\
j_!j^*X \ar[r]^{\epsilon_{_X}}& X\ar[r]^-{f} & i_*i^*X\ar[r]& }\] we
see that $j_!j^*X \stackrel {\epsilon_{_X}}\rightarrow X \stackrel
{\eta_{_X}}\rightarrow i_*i^*X\rightarrow$ is a distinguished
triangle.

\vskip5pt

In order to obtain the second distinguished triangle, we embed the
unit morphism $\zeta_{_X}$ into a distinguished triangle $W\stackrel
w \rightarrow X\stackrel {\zeta_{_X}}\rightarrow
j_*j^*X\rightarrow$. Applying $j^*$ we get a distinguished triangle
$j^*W \stackrel {j^*w}\rightarrow j^*X \stackrel
{j^*\zeta_{_X}}\rightarrow j^*j_*j^*X\rightarrow$. Since
$j^*\zeta_{_X}$ is an isomorphism by Remark 1.2$(2)$, we have $j^*W
= 0$. By ${\rm Im}i_\ast = {\rm Ker}j^*$ we have $W = i_*X'$.
Applying $i^!$ to the distinguished triangle $i_*X' \stackrel
w\rightarrow X\stackrel {\zeta_{_X}}\rightarrow j_*j^*X\rightarrow$
and by $i^!j_* = 0$  we know that $i^!w: \ i^!i_*X' \rightarrow
i^!X$ is an isomorphism. Using the unit isomorphism $X' \rightarrow
i^!i_*X'$, we get a distinguished triangle of the form $i_*i^!X
\stackrel {a}\rightarrow X\stackrel {\zeta_{_X}}\rightarrow
j_*j^*X\rightarrow$ with $i^!a$ an isomorphism. It follows that
${\rm Im}j_* = {\rm Ker}i^!$.

\vskip5pt

Now since ${\rm Im}j_* = {\rm Ker}i^!$ and ${\rm Im}i_\ast = {\rm
Ker}j^*$, it follows that we can replace $(i^*, i_*)$ by $(j^*,
j_*)$, and replace $(j_!, j^*)$ by $(i_*, i^!)$, in the
distinguished triangle $j_!j^*X \stackrel {\epsilon_{_X}}\rightarrow
X \stackrel {\eta_{_X}}\rightarrow i_*i^*X\rightarrow$. In this way
we get the second distinguished triangle $i_*i^!X\stackrel
{\omega_{_X}}\rightarrow X\stackrel {\zeta_{_X}}\rightarrow
j_*j^*X\rightarrow.$ $\s$

\section{\bf \ Upper-symmetric recollements}

\subsection{} Given a
recollement of $\mathcal C$ relative to $\mathcal C'$ and $\mathcal
C''$, one usually can {\bf not} expect a recollement of $\mathcal C$
relative to $\mathcal C''$ and $\mathcal C'$. Inspired by [J] we
define

\begin{defn} \ ([J]) \ A triangulated category recollement
\begin{center}
\begin{picture}(120,45)
\put(0,20){\makebox(3,1){${\mathcal C}'$}}
\put(40,30){\vector(-1,0){30}} \put(10,20){\vector(1,0){30}}
\put(40,10){\vector(-1,0){30}} \put(50,20){\makebox(3,1){${\mathcal
C}$}} \put(90,30){\vector(-1,0){30}} \put(60,20){\vector(1,0){30}}
\put(90,10){\vector(-1,0){30}} \put(100,20){\makebox(3,1){${\mathcal
C}^{\prime\prime}$}} \put(25,35){\makebox(3,1){\scriptsize$i^\ast$}}
\put(25,25){\makebox(3,1){\scriptsize$i_\ast$}}
\put(25,5){\makebox(3,1){\scriptsize$i^!$}}
\put(75,35){\makebox(3,1){\scriptsize$j_!$}}
\put(75,25){\makebox(3,1){\scriptsize$j^\ast$}}
\put(75,5){\makebox(3,1){\scriptsize$j_\ast$}}
\put(240,15){\makebox(25,1)}\put(255,15){\makebox(25,1){$(2.1)$}}
\end{picture}
\end{center}
\noindent of  $\mathcal C$ is upper-symmetric, if there are exact
functors $j_?$ and $i_?$ such that

\begin{center}
\begin{picture}(120,45)
\put(0,20){\makebox(3,1){${\mathcal C}''$}}
\put(40,30){\vector(-1,0){30}} \put(10,20){\vector(1,0){30}}
\put(40,10){\vector(-1,0){30}} \put(50,20){\makebox(3,1){${\mathcal
C}$}} \put(90,30){\vector(-1,0){30}} \put(60,20){\vector(1,0){30}}
\put(90,10){\vector(-1,0){30}} \put(100,20){\makebox(3,1){${\mathcal
C}'$}}\put(25,35){\makebox(3,1){\scriptsize$j^\ast$}}
\put(25,25){\makebox(3,1){\scriptsize$j_\ast$}}
\put(25,5){\makebox(3,1){\scriptsize$j_?$}}
\put(75,35){\makebox(3,1){\scriptsize$i_*$}}
\put(75,25){\makebox(3,1){\scriptsize$i^!$}}
\put(75,5){\makebox(3,1){\scriptsize$i_?$}}
\put(275,15){\makebox(25,1)} \put(255,15){\makebox(25,1){$(2.2)$}}
\end{picture}
\end{center}
\noindent is a recollement; and it is lower-symmetric, if there are
exact functors $j^?$ and $i^?$ such that

\begin{center}
\begin{picture}(120,45)
\put(0,20){\makebox(3,1){${\mathcal C}''$}}
\put(40,30){\vector(-1,0){30}} \put(10,20){\vector(1,0){30}}
\put(40,10){\vector(-1,0){30}} \put(50,20){\makebox(3,1){${\mathcal
C}$}} \put(90,30){\vector(-1,0){30}} \put(60,20){\vector(1,0){30}}
\put(90,10){\vector(-1,0){30}} \put(100,20){\makebox(3,1){${\mathcal
C}'$}}\put(25,35){\makebox(3,1){\scriptsize$j^?$}}
\put(25,25){\makebox(3,1){\scriptsize$j_!$}}
\put(25,5){\makebox(3,1){\scriptsize$j^*$}}
\put(75,35){\makebox(3,1){\scriptsize$i^?$}}
\put(75,25){\makebox(3,1){\scriptsize$i^*$}}
\put(75,5){\makebox(3,1){\scriptsize$i_*$}}
\put(275,15){\makebox(25,1)} \put(255,15){\makebox(25,1){$(2.3)$}}
\end{picture}
\end{center}
\noindent is a recollement. A recollement is symmetric if it is
upper- and lower-symmetric.

\vskip5pt

Similarly, we have a (upper-, lower-) symmetric abelian category
recollement, and note that in abelian situations, all the involved
functors, in particular $j_?$, $i_?$, $j^?$ and $i^?$,  are only
required to be {\bf additive} functors, {\bf not} required to be
exact.
\end{defn}

Let $k$ be a field. P. J$\o$rgensen [J] observed that if a
Hom-finite $k$-linear triangulated category $\mathcal C$ has a Serre
functor, then any recollement of $\mathcal C$ is symmetric: his
proof does not use any triangulated structure of $\mathcal C$ and
hence also works for a Hom-finite $k$-linear abelian category having
a Serre functor. For a similar notion of symmetric recollements of
unbounded derived categories we refer to S. K\"onig [K], and also
Chen-Lin [CL].

\vskip10pt

\subsection{} \  Given two triangulated or abelian category recollements
\begin{center}
\begin{picture}(260,50)
\put(0,20){\makebox(3,1){${\mathcal C}'$}}
\put(40,35){\vector(-1,0){30}} \put(10,20){\vector(1,0){30}}
\put(40,5){\vector(-1,0){30}} \put(50,20){\makebox(3,1){${\mathcal
C}$}} \put(90,35){\vector(-1,0){30}} \put(60,20){\vector(1,0){30}}
\put(90,5){\vector(-1,0){30}} \put(100,20){\makebox(3,1){${\mathcal
C}^{\prime\prime}$}} \put(135,20){\makebox(3,1){\mbox{and}}}
\put(25,40){\makebox(3,1){\scriptsize$i^\ast$}}
\put(25,25){\makebox(3,1){\scriptsize$i_\ast$}}
\put(25,10){\makebox(3,1){\scriptsize$i^!$}}
\put(75,40){\makebox(3,1){\scriptsize$j_!$}}
\put(75,25){\makebox(3,1){\scriptsize$j^\ast$}}
\put(75,10){\makebox(3,1){\scriptsize$j_\ast$}}

\put(165,20){\makebox(3,1){${\mathcal C}'$}}
\put(205,35){\vector(-1,0){30}} \put(175,20){\vector(1,0){30}}
\put(205,5){\vector(-1,0){30}} \put(215,20){\makebox(3,1){${\mathcal
D}$}} \put(255,35){\vector(-1,0){30}} \put(225,20){\vector(1,0){30}}
\put(255,5){\vector(-1,0){30}} \put(265,20){\makebox(3,1){${\mathcal
C}^{\prime\prime}$}}
\put(190,40){\makebox(3,1){\scriptsize$i_D^\ast$}}
\put(190,25){\makebox(3,1){\scriptsize$i^D_\ast$}}
\put(190,10){\makebox(3,1){\scriptsize$i_D^!$}}
\put(240,40){\makebox(3,1){\scriptsize$j^D_!$}}
\put(240,25){\makebox(3,1){\scriptsize$j_D^\ast$}}
\put(240,10){\makebox(3,1){\scriptsize$j^D_\ast$}}
\end{picture}
\end{center}
\noindent if there is an exact functor $f: {\mathcal C}\rightarrow
{\mathcal D}$ such that there are natural isomorphisms
$$i^\ast \approx i^\ast_D f, \ \  f i_\ast\approx
i^D_\ast, \ \ i^!\approx i^!_D f, \ \ fj_!\approx j^D_!, \ \
j^\ast\approx j^\ast_D f, \ \ fj_\ast\approx j^D_\ast,$$ then we
call $f$ {\it a comparison functor}. Two (triangulated or abelian
category) recollements are {\it equivalent} if there is a comparison
functor $f$ which is an equivalence of categories. According to
Parshall-Scott [PS, Theorem 2.5], a comparison functor between
triangulated category recollements is an equivalence of categories.
However, Franjou-Pirashvili [FV] pointed out that this is not
necessarily the case for abelian category recollements.

\vskip10pt

\subsection{} In this subsection we only consider triangulated category
recollements.  If $(2.1)$ is an upper-symmetric recollement, then we
call $(2.2)$ {\it a upper-symmetric version} of $(2.1)$; and if
$(2.1)$ is an lower-symmetric recollement, then we call $(2.3)$ {\it
a lower-symmetric version} of $(2.1)$.

\vskip10pt

\begin{lem} \ $(1)$ \ Any two upper-symmetric versions of a upper-symmetric recollement
are equivalent.

\vskip5pt

$(1')$ \ \ Any two lower-symmetric versions of a lower-symmetric
recollement are equivalent.

\vskip5pt

$(2)$ \ Equivalent upper-symmetric recollements have equivalent
upper-symmetric versions.

\vskip5pt

$(2')$ \  Equivalent lower-symmetric recollements have equivalent
lower-symmetric versions.
\end{lem} \noindent{\bf Proof.} \ $(1)$ \ Let
$(2.2)$ and \begin{center}
\begin{picture}(120,45)
\put(0,20){\makebox(3,1){${\mathcal C}''$}}
\put(40,30){\vector(-1,0){30}} \put(10,20){\vector(1,0){30}}
\put(40,10){\vector(-1,0){30}} \put(50,20){\makebox(3,1){${\mathcal
C}$}} \put(90,30){\vector(-1,0){30}} \put(60,20){\vector(1,0){30}}
\put(90,10){\vector(-1,0){30}} \put(100,20){\makebox(3,1){${\mathcal
C}'$}}\put(25,35){\makebox(3,1){\scriptsize$j^\ast$}}
\put(25,25){\makebox(3,1){\scriptsize$j_\ast$}}
\put(25,5){\makebox(3,1){\scriptsize$j_{??}$}}
\put(75,35){\makebox(3,1){\scriptsize$i_*$}}
\put(75,25){\makebox(3,1){\scriptsize$i^!$}}
\put(75,5){\makebox(3,1){\scriptsize$i_{??}$}}
\put(275,15){\makebox(25,1)}\put(255,15){\makebox(25,1){$(2.4)$}}
\end{picture}
\end{center}
be two upper-symmetric versions of a upper-symmetric recollement
$(2.1)$. Then $j_{??}i_? = 0.$ In fact, for $Y\in\mathcal C'$ we
have
$$\Hom_{\mathcal C''}(j_{??}i_?Y, j_{??}i_?Y) \cong\Hom_{\mathcal
C}(j_*j_{??}i_?Y, i_?Y) \cong \Hom_{\mathcal C''}(i^!j_*j_{??}i_?Y,
Y) = 0.$$ For $X\in \mathcal C$, by $(2.2)$ and ${\rm (R4)}$ we have
distinguished triangle $j_*j_?X \stackrel {\epsilon_{_X}}\rightarrow
X \stackrel {\eta_{_X}}\rightarrow i_?i^!X \rightarrow.$ Applying
exact functor $j_{??}$ and using the unit ${\rm Id}_{\mathcal
C''}\rightarrow j_{??}j_*$, we have
$$j_{??}X \cong j_{??}j_*j_?X\cong j_?X,$$
which means that $j_{??}$ is naturally isomorphic to $j_?$.
Similarly one can prove that $i_{??}$ is naturally isomorphic to
$i_?$. Thus ${\rm Id}_{\mathcal C}$ is an equivalence between
$(2.2)$ and $(2.4)$. This proves $(1)$.

$(1')$ can be similarly proved.

\vskip5pt

$(2)$ \ Given two equivalent upper-symmetric recollements
\begin{center}
\begin{picture}(260,50)
\put(0,20){\makebox(3,1){${\mathcal C}'$}}
\put(40,35){\vector(-1,0){30}} \put(10,20){\vector(1,0){30}}
\put(40,5){\vector(-1,0){30}} \put(50,20){\makebox(3,1){${\mathcal
C}$}} \put(90,35){\vector(-1,0){30}} \put(60,20){\vector(1,0){30}}
\put(90,5){\vector(-1,0){30}} \put(100,20){\makebox(3,1){${\mathcal
C}^{\prime\prime}$}} \put(135,20){\makebox(3,1){\mbox{and}}}
\put(25,40){\makebox(3,1){\scriptsize$i^\ast$}}
\put(25,25){\makebox(3,1){\scriptsize$i_\ast$}}
\put(25,10){\makebox(3,1){\scriptsize$i^!$}}
\put(75,40){\makebox(3,1){\scriptsize$j_!$}}
\put(75,25){\makebox(3,1){\scriptsize$j^\ast$}}
\put(75,10){\makebox(3,1){\scriptsize$j_\ast$}}

\put(165,20){\makebox(3,1){${\mathcal C}'$}}
\put(205,35){\vector(-1,0){30}}
\put(175,20){\vector(1,0){30}}
\put(205,5){\vector(-1,0){30}}
\put(215,20){\makebox(3,1){${\mathcal
D}$}}
\put(255,35){\vector(-1,0){30}}
\put(225,20){\vector(1,0){30}}
\put(255,5){\vector(-1,0){30}}
\put(265,20){\makebox(3,1){${\mathcal
C}^{\prime\prime}$}}
\put(190,40){\makebox(3,1){\scriptsize$i_D^\ast$}}
\put(190,25){\makebox(3,1){\scriptsize$i^D_\ast$}}
\put(190,10){\makebox(3,1){\scriptsize$i_D^!$}}
\put(240,40){\makebox(3,1){\scriptsize$j^D_!$}}
\put(240,25){\makebox(3,1){\scriptsize$j_D^\ast$}}
\put(240,10){\makebox(3,1){\scriptsize$j^D_\ast$}}
\end{picture}
\end{center}

\noindent with comparison functor $f$, let $(2.2)$ as an
upper-symmetric version of the first recollement.  By Lemma 1.3 we
know that
\begin{center}
\begin{picture}(120,50)
\put(0,20){\makebox(3,1){${\mathcal C}''$}}
\put(40,35){\vector(-1,0){30}} \put(10,20){\vector(1,0){30}}
\put(40,5){\vector(-1,0){30}} \put(50,20){\makebox(3,1){${\mathcal
D}$}} \put(90,35){\vector(-1,0){30}} \put(60,20){\vector(1,0){30}}
\put(90,5){\vector(-1,0){30}} \put(100,20){\makebox(3,1){${\mathcal
C'}$} \ $_\cdot$} \put(25,40){\makebox(3,1){\scriptsize$j_D^\ast$}}
\put(25,25){\makebox(3,1){\scriptsize$j^D_\ast$}}
\put(25,10){\makebox(3,1){\scriptsize$j_{?}f^{-1}$}}
\put(75,40){\makebox(3,1){\scriptsize$i^D_*$}}
\put(75,25){\makebox(3,1){\scriptsize$i_D^!$}}
\put(75,10){\makebox(3,1){\scriptsize$fi_{?}$}}
\put(275,15){\makebox(25,1)}\put(255,15){\makebox(25,1){$(2.5)$}}
\end{picture}
\end{center}
is a triangulated category recollement, and that $f$ is an
equivalence between $(2.2)$ and $(2.5)$. Note that $(2.5)$ is an
upper-symmetric version of the second given upper-symmetric
recollement, and hence the assertion follows from $(1)$.

$(2')$ can be similarly proved. $\s$

\vskip10pt

Let $\mathcal C', \mathcal C, \mathcal C''$ be triangulated
categories. Denote by $USR(\mathcal C', \mathcal C, \mathcal C'')$
the class of equivalence classes of the upper-symmetric recollements
of triangulated category $\mathcal C$ relative to $\mathcal C'$ and
$\mathcal C''$; and denote by $LSR(\mathcal C', \mathcal C, \mathcal
C'')$ the class of the lower-symmetric recollements of triangulated
category $\mathcal C$ relative to $\mathcal C''$ and $\mathcal C'$.

\vskip10pt

\begin{thm} \ There is a one-one correspondence between
$USR(\mathcal C', \mathcal C, \mathcal C'')$ and $LSR(\mathcal C'',
\mathcal C, \mathcal C')$.\end{thm} \noindent{\bf Proof.} \ Given an
upper-symmetric recollement $(2.1)$, observe that an upper-symmetric
version $(2.2)$ of $(2.1)$ is lower-symmetric: in fact, $(2.1)$
could be a lower-symmetric version of $(2.2)$.  Similarly, a
lower-symmetric recollement could be a upper-symmetric version of a
lower-symmetric version of itself. Thus by Lemma 2.2 we get a
one-one correspondence between $USR(\mathcal C', \mathcal C,
\mathcal C'')$ and $LSR(\mathcal C'', \mathcal C, \mathcal C')$.
$\s$

\subsection{} \ We consider Artin algebras over a fixed commutative artinian
ring, and finitely generated modules. Let $A$ and $B$ be Artin
algebras, and $M$ an $A$-$B$-bimodule. Then $\m =
\left(\begin{smallmatrix}
A&M\\
0&B
\end{smallmatrix}\right)$ is an Artin algebra with multiplication given by the one of matrices.
Denoted by $A$-mod the category of finitely generated left
$A$-modules. A left $\m$-module is identified with a triple
$\left(\begin{smallmatrix}
X\\
Y
\end{smallmatrix}\right)_\phi$, or simply  $\left(\begin{smallmatrix}
X\\
Y
\end{smallmatrix}\right)$ if $\phi$ is clear, where $X\in A$-mod, $Y\in B$-mod,
and $\phi: M\otimes_{B} Y\rightarrow X$ is an
$A$-map. A $\m$-map $\left(\begin{smallmatrix} X \\
   Y
 \end{smallmatrix}\right)_{\phi}\rightarrow \left(\begin{smallmatrix}
   X'  \\
   Y'
 \end{smallmatrix}\right)_{\phi'}$ is identified with a pair $\left(\begin{smallmatrix}
   f  \\
   g
 \end{smallmatrix}\right)$, where $f\in {\Hom}_A (X, X'),$ $g\in\Hom_B (Y, Y')$,
such that $\phi'({\rm Id}\otimes g) = f\phi$. The indecomposable
projective $\m$-modules are exactly $\left(\begin{smallmatrix}
   P \\
   0
\end{smallmatrix}\right)
$ and $\left( \begin{smallmatrix}
   M\otimes_B Q  \\
   Q
\end{smallmatrix}\right)_{\rm id}$,
where $P$ runs over indecomposable projective $A$-modules, and $Q$
runs over indecomposable projective $B$-modules. See [ARS], p.73.

For any $A$-module $X$ and $B$-module $Y$, denote by
$\alpha_{_{X,Y}}$ the adjoint isomorphism
$$\alpha_{_{X,Y}}: \ \Hom_A(M\otimes_B Y, X)\longrightarrow \Hom_B(Y, \Hom_A(M,
X))$$ given by $$\alpha_{_{X,Y}}(\phi)(y)(m) = \phi(m\otimes y), \
\forall \ \phi\in\Hom_A(M\otimes_B Y, X), \ y\in Y, \ m\in M.$$ Put
$\psi_{_X}$ to be $\alpha^{-1}_{_{X, \Hom(M, X)}}({\rm
Id}_{_{\Hom(M, X)}})$. Thus $\psi_{_X}: M\otimes_B\Hom_A(M,
X)\rightarrow X$ is given by $m\otimes f\mapsto f(m)$.

\vskip10pt

\begin{thm} \ Let $A$ and $B$ be Artin algebras,
$_AM_B$ an $A$-$B$-bimodule,  and $\m = \left(\begin{smallmatrix}
A&M\\
0&B
\end{smallmatrix}\right)$. Then we have an upper-symmetric (but non
lower-symmetric) abelian category recollement

\begin{center}
\begin{picture}(120,45)
\put(0,20){\makebox(-16,1){$A\mbox{-}{\rm mod}$}}
\put(40,30){\vector(-1,0){30}} \put(10,20){\vector(1,0){30}}
\put(40,10){\vector(-1,0){30}} \put(50,20) {\makebox(16,1)
{$\m\mbox{-}{\rm mod}$}} \put(105,30){\vector(-1,0){30}}
\put(75,20){\vector(1,0){30}} \put(105,10){\vector(-1,0){30}}
\put(115,20){\makebox(18, 0.5) {$B\mbox{-}{\rm mod}$}}
\put(25,35){\makebox(3,1){\scriptsize$i^\ast$}}
\put(25,25){\makebox(3,1){\scriptsize$i_\ast$}}
\put(25,5){\makebox(3,1){\scriptsize$i^!$}}
\put(90,35){\makebox(3,1){\scriptsize$j_!$}}
\put(90,25){\makebox(3,1){\scriptsize$j^\ast$}}
\put(90,5){\makebox(3,1){\scriptsize$j_\ast$}}
\put(250,15){\makebox(25,1){$(2.6)$}}\end{picture}
\end{center}
where

\vskip5pt$i^*$ is given by $\left(\begin{smallmatrix}
X\\
Y
\end{smallmatrix}\right)_\phi\mapsto {\rm \Cok}\phi$;
\ \ $i_*$ is given by $X\mapsto \left(\begin{smallmatrix}
X\\
0
\end{smallmatrix}\right)$; \ \ $i^!$ is given by $\left(\begin{smallmatrix}
X\\
Y
\end{smallmatrix}\right)_\phi\mapsto X$;

\vskip5pt  $j_!$ is given by $Y\mapsto \left(\begin{smallmatrix}
M\otimes Y\\
Y
\end{smallmatrix}\right)_{\rm Id}$; \ \ $j^*$ is given by $\left(\begin{smallmatrix}
X\\
Y
\end{smallmatrix}\right)_{\phi}\mapsto Y$; \ \ $j_*$ is given by $Y\mapsto \left(\begin{smallmatrix}
0\\
Y
\end{smallmatrix}\right);$

\vskip5pt$j_?$ is given by $\left(\begin{smallmatrix}
X\\
Y
\end{smallmatrix}\right)_\phi\mapsto \Ker\alpha_{_{X,Y}}(\phi);$
\ and  \ $i_?$ is given by $X\mapsto \left(\begin{smallmatrix}
X\\
\Hom_A(M, X)
\end{smallmatrix}\right)_{\psi_{_X}}.$
\end{thm}
\noindent{\bf Proof.} \ By construction $i_*, \ j_!$ and $j_*$ are
fully faithful;  ${\rm Im}i_\ast = {\rm Ker}j^*$, and ${\rm
Im}j_\ast = {\rm Ker}i^!$. For  $\left(\begin{smallmatrix}
X\\
Y
\end{smallmatrix}\right)_\phi\in \m$-mod,
$X'\in A$-mod, and $Y'\in B$-mod, we have the following isomorphisms
of abelian groups, which are natural in both positions
$$\Hom_A (\Cok\phi, X') \cong \Hom_\m(\left(\begin{smallmatrix}
X\\
Y
\end{smallmatrix}\right)_\phi, \left(\begin{smallmatrix}
X'\\
0
\end{smallmatrix}\right))
\eqno(2.7)$$ given by $f\mapsto \left(\begin{smallmatrix}
f\pi\\
0
\end{smallmatrix}\right)$, where $\pi: X\rightarrow \Cok\phi$ is
the canonical $A$-map;
$$\Hom_\m(\left(\begin{smallmatrix}
X'\\
0
\end{smallmatrix}\right), \left(\begin{smallmatrix}
X\\
Y
\end{smallmatrix}\right)_\phi)
\cong \Hom_A (X', X);  \eqno(2.8)$$
$$\Hom_\m(\left(\begin{smallmatrix}
M\otimes Y'\\
Y'
\end{smallmatrix}\right)_{\rm Id}, \left(\begin{smallmatrix}
X\\
Y
\end{smallmatrix}\right)_{\phi}) \cong \Hom_B (Y', Y) \eqno(2.9)$$
given by $\left(\begin{smallmatrix}
\phi({\rm Id}\otimes g)\\
g
\end{smallmatrix}\right)\mapsto g$; \ and
$$\Hom_B (Y,
Y') \cong \Hom_\m(\left(\begin{smallmatrix}
X\\
Y
\end{smallmatrix}\right)_{\phi}, \left(\begin{smallmatrix}
0\\
Y'
\end{smallmatrix}\right)).$$
Thus $(i^\ast, i_\ast)$,  \ $(i_\ast, i^!)$, \ $(j_!, j^\ast)$, and
$(j^\ast, j_\ast)$ are adjoint pairs, and hence $(2.6)$ is a
recollement. It is not lower-symmetric since ${\rm Im}j_! \ne {\rm
Ker}i^*$.

In order to see that it is upper-symmetric, it remains to prove that
$(j_*, j_?)$ and $(i^!, i_?)$ are adjoint pairs, and that $i_?$ is
fully faithful. For $g\in\Hom_B(Y, Y')$ and
$\left(\begin{smallmatrix}
X'\\
Y'
\end{smallmatrix}\right)_{\phi'}\in\m$-mod, we have
\begin{align*} & \left(\begin{smallmatrix}
0\\
g
\end{smallmatrix}\right)\in \Hom_\m(\left(\begin{smallmatrix}
0\\
Y
\end{smallmatrix}\right), \left(\begin{smallmatrix}
X'\\
Y'
\end{smallmatrix}\right)_{\phi'})  \Longleftrightarrow \phi'({\rm Id}\otimes g) =
0 \Longleftrightarrow  \phi'(m\otimes g(y)) = 0, \ \forall \ y\in Y,
\ \forall \ m\in M
\\ &\Longleftrightarrow  \alpha_{_{X', Y'}}(\phi')(g(y)) = 0, \ \forall \ y\in Y\Longleftrightarrow  g(Y)\subseteq \Ker\alpha_{_{X', Y'}}(\phi')
\Longleftrightarrow  g\in\Hom_B(Y, \Ker\alpha_{_{X', Y'}}(\phi')).
\end{align*}
It follows that $\left(\begin{smallmatrix}
0\\
g
\end{smallmatrix}\right)\mapsto g$ gives an isomorphism
$\Hom_\m(\left(\begin{smallmatrix}
0\\
Y
\end{smallmatrix}\right), \left(\begin{smallmatrix}
X'\\
Y'
\end{smallmatrix}\right)_{\phi'})
\rightarrow \Hom_B(Y, \Ker\alpha_{_{X', Y'}}(\phi'))$ of abelian
groups, which is natural in both positions, i.e., $(j_*, j_?)$ is an
adjoint pair. Let $\left(\begin{smallmatrix}
f\\
g
\end{smallmatrix}\right)\in \Hom_\m(\left(\begin{smallmatrix}
X\\
Y
\end{smallmatrix}\right)_\phi, \left(\begin{smallmatrix}
X'\\
\Hom_A(M, X')
\end{smallmatrix}\right)_{\psi_{_{X'}}})$. By  $\psi_{{_{X'}}}({\rm id}\otimes g) = f\phi$
we have
\begin{align*}\alpha_{_{X', Y}}(f\phi)(y)(m) &=  f\phi(m\otimes y) = \psi_{_{X'}}({\rm Id}\otimes g)(m\otimes y)
\\ & =\psi_{_{X'}}(m\otimes g(y)) = g(y)(m), \ \forall \ y\in Y, \ \forall \ m\in M, \end{align*}
which means $g = \alpha_{_{X', Y}}(f\phi).$ Thus $f\mapsto
\left(\begin{smallmatrix}
f\\
\alpha_{_{X', Y}}(f\phi)
\end{smallmatrix}\right)$ gives an isomorphism
$$\Hom_A(X, X')\longrightarrow\Hom_\m(\left(\begin{smallmatrix}
X\\
Y
\end{smallmatrix}\right)_\phi, \left(\begin{smallmatrix}
X'\\
\Hom_A(M, X')
\end{smallmatrix}\right)_{\psi_{_{X'}}})
$$ of abelian groups, which is natural in both
positions, i.e., $(i^!, i_?)$ is an adjoint pair. Since

\noindent$\alpha_{_{X', \Hom(M, X)}}(f\psi_{_X}) = \Hom_A(M, f)$,
this isomorphism also shows that $i_?$ is fully faithful. This
completes the proof. $\s$

\vskip10pt

By Theorem 2.4 we have

\vskip10pt

\begin{cor} \ Let $A$ be a Gorenstein algebra, and $T_2(A)
= \left(\begin{smallmatrix}
A&A\\
0&A
\end{smallmatrix}\right)$.
Then we have an upper-symmetric (but non lower-symmetric) abelian
category recollement

\begin{center}
\begin{picture}(120,45)
\put(0,20){\makebox(-16,1){$A\mbox{-}{\rm mod}$}}
\put(40,30){\vector(-1,0){30}} \put(10,20){\vector(1,0){30}}
\put(40,10){\vector(-1,0){30}} \put(60,20) {\makebox(16,1)
{$T_2(A)\mbox{-}{\rm mod}$}} \put(125,30){\vector(-1,0){30}}
\put(95,20){\vector(1,0){30}} \put(125,10){\vector(-1,0){30}}
\put(135,20){\makebox(18, 0.5) {$A\mbox{-}{\rm mod.}$}}
\end{picture}
\end{center}
\end{cor}

\begin{rem} \ As we see from
$(2.6)$ and its upper symmetric version, in an abelian category
recollement, the following statement may {\bf not} true:

\vskip5pt

$(1)$  \ ${\rm Im}j_! = {\rm Ker}i^*$; \ ${\rm Im}j_\ast = {\rm
Ker}i^!$;

\vskip5pt

$(2)$ \ The counits and units give rise to exact sequences of
natural transformations: $$0\longrightarrow j_!j^* \longrightarrow
{\rm Id}_{\mathcal C}\longrightarrow i_*i^*\longrightarrow 0 \ \ \
\mbox{and} \ \ \ 0\longrightarrow i_*i^! \longrightarrow {\rm
Id}_{\mathcal C}\longrightarrow j_*j^*\longrightarrow 0.$$

\vskip5pt

$(3)$ \ $i^!j_! = 0$; and  $i^\ast j_\ast = 0$.

\vskip5pt

In triangulated situations,   $(1)$ and the corresponding version of
$(2)$ {\bf always} hold; but $(3)$ is also {\bf not} true in
general.
\end{rem}

\section{\bf \ Symmetric recollements induced by Gorenstein-projective modules}

\subsection{} \ Let $A$  be an Artin algebra. An $A$-module $G$ is {\it Gorenstein-projective}, if there
is an exact sequence $\cdots \rightarrow P^{-1}\rightarrow P^{0}
\stackrel{d^0}{\rightarrow} P^{1}\rightarrow \cdots$ of projective
$A$-modules, which stays exact under ${\rm Hom}_A(-, A)$, and such
that $G\cong \operatorname{Ker}d^0$. Let $A\mbox{-}\mathcal Gproj$
be the full subcategory of $A$-mod consisting  of the
Gorenstein-projective modules. Then $A\mbox{-}\mathcal Gproj
\subseteq \ ^\perp A$, where $^\perp A = \{X\in A\mbox{-}{\rm mod} \
| \ {\rm Ext}^i_A(X, A) = 0,  \forall \ i\ge 1\}$; and $\Hom_A(-, \
_AA)$ induces a duality $A$-${\mathcal Gproj}\cong
A^{op}$-${\mathcal Gproj}$ with a quasi-inverse $\Hom_{A}(-, A_A)$
([B], Proposition 3.4). An important feature is that
$A\mbox{-}\mathcal Gproj$ is a Frobenius category with
projective-injective objects being projective $A$-modules, and hence
the stable category $\underline {A\mbox{-}\mathcal Gproj}$ modulo
projective $A$-modules is a triangulated category ([Hap]).

\vskip10pt

An Artin algebra $A$ is {\it Gorenstein}, if ${\rm inj.dim}\ _AA<
\infty$ and ${\rm inj.dim} \ A_A < \infty$. We have the following
well-known fact (E. Enochs - O. Jenda [EJ], Corollary 11.5.3).

\vskip10pt

\begin {lem}\ Let $A$ be a Gorenstein algebra. Then

\vskip5pt

$(1)$ \ If $P^\bullet$ is an exact sequence of projective left
(resp. right) $A$-modules, then ${\rm Hom}_A(P^\bullet, \ A)$ is
again an exact sequence of projective right (resp. left)
$A$-modules.

\vskip5pt

$(2)$ \ A module $G$ is Gorenstein-projective if and only if there
is an exact sequence $0\rightarrow G\rightarrow P^0\rightarrow
P^1\rightarrow \cdots$ with each $P^i$ projective.

\vskip5pt

$(3)$ \ $A\mbox{-}\mathcal Gproj = \ ^\perp A$.
\end{lem}
\noindent{\bf Proof.} \ For convenience we include an alternating
proof.

\vskip5pt

$(1)$ \ Let $0\rightarrow K \rightarrow I_0\rightarrow
I_1\rightarrow 0$ be an exact sequence with $I_0, \ I_1$ injective
modules. Then  $0\rightarrow {\rm Hom}_A(P^\bullet, \ K) \rightarrow
{\rm Hom}_A(P^\bullet, \ I_0) \rightarrow {\rm Hom}_A(P^\bullet, \
I_1) \rightarrow 0$ is an exact sequence of complexes. Since ${\rm
Hom}_A(P^\bullet, \ I_i) \ (i=0, 1)$ are exact, it follows that
${\rm Hom}_A(P^\bullet, \ K)$ is exact. Repeating this process, by
${\rm inj.dim}\ _AA< \infty$ we deduce that ${\rm Hom}_A(P^\bullet,
\ A)$ is exact.

\vskip5pt

$(2)$ \ This follows from definition and $(1)$.

\vskip10pt

$(3)$ \ Let $G\in \ ^\perp A$. Applying $\Hom_A(-, A)$ to a
projective resolution of $G$ we get is an exact sequence. By $(2)$
this means that $\Hom_A(G, A)$ is a Gorenstein-projective right
$A$-module, and hence $G$ is Gorenstein-projective by the duality
$\Hom_A(-, \ _AA): \ A\mbox{-}{\mathcal Gproj}\cong
A^{op}\mbox{-}{\mathcal Gproj}$. $\s$

\vskip10pt

We need the following description of Gorenstein-projective
$\m$-modules.

\vskip10pt

\begin{prop} \ Let $A$ and $B$ be Gorenstein algebras,
$M$ an $A$-$B$-bimodule such that $_AM$ and $M_B$ are projective,
and $\m = \left(\begin{smallmatrix}
A&M\\
0&B
\end{smallmatrix}\right)$. Then
$\left(\begin{smallmatrix}X\\Y\end{smallmatrix}\right)_\phi$ is a
Gorenstein-projective $\Lambda$-module if and only if $\phi:
M\otimes Y\rightarrow X$ is monic, $X$ and $\Cok\phi$ are
Gorenstein-projective $A$-modules, and $Y$ is a
Gorenstein-projective $B$-module. In this case $M\otimes Y$ is a
Gorenstein-projective $A$-module. \
\end{prop} \noindent{\bf Proof.} \ If
$\left(\begin{smallmatrix}X\\Y\end{smallmatrix}\right)_\phi$ is a
Gorenstein-projective $\Lambda$-module, then there is an exact
sequence
$$0\longrightarrow \left(\begin{smallmatrix}X\\Y\end{smallmatrix}\right)_\phi\longrightarrow
\left(\begin{smallmatrix}P_0\oplus (M\otimes
Q_0)\\Q_0\end{smallmatrix}\right)_{\binom{0}{\rm Id}}\longrightarrow
\left(\begin{smallmatrix}P_1\oplus (M\otimes
Q_1)\\Q_1\end{smallmatrix}\right)_{\binom{0}{\rm
Id}}\longrightarrow\cdots\eqno(3.1)$$ where $P_i$ and $Q_i$ are
respectively projective $A$- and $B$-modules,  $i\ge 0$, i.e., we
have exact sequences
$$0\longrightarrow X \longrightarrow P_0\oplus (M\otimes
Q_0)\longrightarrow P_1\oplus (M\otimes Q_1)\longrightarrow
\cdots\eqno(3.2)$$ and
$$0\longrightarrow Y \longrightarrow Q_0\longrightarrow Q_1\longrightarrow \cdots,\eqno(3.3)$$
such that the following diagram commutes
\[\xymatrix {\ \ \ \ \ \ \ \
0 \ar[r] & M\otimes_BY\ar[d]^-{\phi}\ar[r]& M\otimes_B
Q_0\ar[d]_{\binom{0}{\rm Id}}\ar[r]& M\otimes_B Q_1
\ar[d]^-{\binom{0}{\rm Id}}\ar[r] & \cdots \ \ \ \ \ \ \ \ \ \ \ \ \
\ \ \ \ \ \ \ \ \ \ \    \\ \ \ \ \ \ \ \ \ 0 \ar[r] & X\ar[r] &
P_0\oplus (M\otimes Q_0)\ar[r]& P_1\oplus (M\otimes Q_1)\ar[r] &
\cdots. \ \ \ \ \ \ \ \ \ \ \ \ \ \  \ \ \   (3.4)}\] By Lemma
3.1$(2)$ $Y$ is Gorenstein-projective. Since $_AM$ and $_BQ_i$ are
projective, it follows that $M\otimes Q_i$ are projective
$A$-modules, and hence $X$ is Gorenstein-projective by Lemma
3.1$(2)$. Since $M_B$ is projective, by $(3.3)$ the upper row of
$(3.4)$ is exact, and hence $M\otimes Y$ is Gorenstein-projective
and $\phi$ is monic. By $(3.4)$ we get exact sequence $0\rightarrow
\Cok\phi\rightarrow P_0\rightarrow P_1\rightarrow \cdots,$ thus
$\Cok\phi$ is Gorenstein-projective by Lemma 3.1$(2)$.

\vskip5pt

Conversely, we have exact sequence $(3.3)$ with $Q_i$ being
projective $B$-modules. Since $M_B$ is projective and $\Cok\phi$ is
Gorenstein-projective, we get the following exact sequences
\begin{align*} & 0 \longrightarrow M\otimes Y\longrightarrow
M\otimes  Q_0\longrightarrow M\otimes Q_1 \longrightarrow \cdots
\\ & 0 \longrightarrow \Cok\phi\longrightarrow \ \ \ P_0 \ \ \ \ \longrightarrow \ \ \   P_1 \ \ \ \ \longrightarrow  \cdots
\end{align*}
with $P_i$ projective. Since $M\otimes Q_i  \ (i\ge 0)$ are
projective $A$-modules and projective $A$-modules are injective
objects in $A$-$\mathcal Gproj$, it follows from the exact sequence
$0\rightarrow M\otimes Y\rightarrow X \rightarrow\Cok\phi\rightarrow
0$ and a version of Horseshoe Lemma that there is an exact sequence
$(3.2)$ such that the diagram $(3.4)$ commutes. This means that
$(3.1)$ is exact. Since $\m$ is also Gorenstein (see e.g. [C],
Theorem 3.3), it follows from Lemma 3.1$(2)$ that
$\left(\begin{smallmatrix}X\\Y\end{smallmatrix}\right)_\phi$ is a
Gorenstein-projective $\Lambda$-module. $\s$

\subsection{} \ The main result of this section is as follows.

\begin{thm} \ Let $A$ and $B$ be Gorenstein algebras,
$M$ an $A$-$B$-bimodule such that $_AM$ and $M_B$ are projective,
and $\m = \left(\begin{smallmatrix}
A&M\\
0&B
\end{smallmatrix}\right)$. Then we have  a triangulated category recollement

\begin{center}
\begin{picture}(120,45)
\put(0,20){\makebox(-22,1){$\underline {A\mbox{-}\mathcal Gproj}$}}
\put(40,30){\vector(-1,0){30}} \put(10,20){\vector(1,0){30}}
\put(40,10){\vector(-1,0){30}}
\put(50,20){\makebox(25,0.8){$\underline {\m\mbox{-}\mathcal
Gproj}$}} \put(115,30){\vector(-1,0){30}}
\put(85,20){\vector(1,0){30}} \put(115,10){\vector(-1,0){30}}
\put(125,20){\makebox(27,0.5){$\underline {B\mbox{-}\mathcal
Gproj}$} \  .}\put(25,35){\makebox(3,1){\scriptsize$i^\ast$}}
\put(25,25){\makebox(3,1){\scriptsize$i_\ast$}}
\put(25,5){\makebox(3,1){\scriptsize$i^!$}}
\put(100,35){\makebox(3,1){\scriptsize$j_!$}}
\put(100,25){\makebox(3,1){\scriptsize$j^\ast$}}
\put(100,5){\makebox(3,1){\scriptsize$j_\ast$}}
\end{picture}
\end{center}
Moreover, if $A$ and $B$ are in additional finite-dimensional
algebras over a field, then it is a symmetric recollement.
\end{thm}

\subsection{} \ Before giving a proof,  we construct all the functors in Theorem 3.3.
If a $\m$-map $\left(\begin{smallmatrix}
X\\
Y
\end{smallmatrix}\right)_\phi\rightarrow \left(\begin{smallmatrix}
X'\\
Y'
\end{smallmatrix}\right)_{\phi'}$ factors through a projective
$\m$-module $\left(\begin{smallmatrix}P\\
0
\end{smallmatrix}\right)\oplus \left(\begin{smallmatrix}
M\otimes Q\\
Q
\end{smallmatrix}\right)$, then it is easy to see that the induced $A$-map ${\rm \Cok}\phi\rightarrow {\rm
\Cok}\phi'$ factors through $P$. By Proposition 3.2 this implies
that the functor ${\m\mbox{-}\mathcal Gproj}\rightarrow \underline{
A\mbox{-}\mathcal Gproj}$ given by $\left(\begin{smallmatrix}
X\\
Y
\end{smallmatrix}\right)_\phi\mapsto {\rm \Cok}\phi$ induces a
functor $i^*: \ \underline {\m\mbox{-}\mathcal Gproj}\rightarrow
\underline{ A\mbox{-}\mathcal Gproj}$.

\vskip5pt

By Proposition 3.2 there is a unique functor $i_*: \underline{
A\mbox{-}\mathcal Gproj}\rightarrow \underline {\m\mbox{-}\mathcal
Gproj}$ given by $X\mapsto \left(\begin{smallmatrix}
X\\
0
\end{smallmatrix}\right)$, which is fully faithful.

\vskip5pt

If a $\m$-map $\left(\begin{smallmatrix}
f\\
g
\end{smallmatrix}\right): \left(\begin{smallmatrix}
X\\
Y
\end{smallmatrix}\right)_\phi\rightarrow \left(\begin{smallmatrix}
X'\\
Y'
\end{smallmatrix}\right)_{\phi'}$ factors through a projective
$\m$-module $\left(\begin{smallmatrix}P\\
0
\end{smallmatrix}\right)\oplus \left(\begin{smallmatrix}
M\otimes Q\\
Q
\end{smallmatrix}\right)$, then $f: \ X\rightarrow X'$ factors through a projective
$A$-module $P\oplus (M\otimes Q)$. By Proposition 3.2 this implies
that there is a unique functor $i^!:  \underline {\m\mbox{-}\mathcal
Gproj}\rightarrow \underline{ A\mbox{-}\mathcal Gproj}$ given by
$\left(\begin{smallmatrix}
X\\
Y
\end{smallmatrix}\right)_\phi\mapsto X$.

\vskip5pt

By Proposition 3.2 there is a unique functor $j^*: \ \underline{
\m\mbox{-}\mathcal Gproj}\rightarrow \underline {B\mbox{-}\mathcal
Gproj}$ given by $\left(\begin{smallmatrix}
X\\
Y
\end{smallmatrix}\right)_{\phi}\mapsto Y$.

\vskip5pt

Let $_BY$ be a Gorenstein-projective module. Since $M_B$ is
projective, by Lemma 3.1$(2)$ $M\otimes Y$ is a
Gorenstein-projective $A$-module.  By Proposition 3.2 there is a
unique functor $j_!: \ \underline{ B\mbox{-}\mathcal
Gproj}\rightarrow \underline {\m\mbox{-}\mathcal Gproj}$ given by
$Y\mapsto \left(\begin{smallmatrix}
M\otimes Y\\
Y
\end{smallmatrix}\right)_{\rm Id}$, which is fully faithful.

\vskip10pt

\begin{lem} \ Let $A$, $B$, $M$, and
$\m$ be as in Theorem 3.3. Then there exists a unique fully faithful
functor $j_*: \ \underline{ B\mbox{-}\mathcal Gproj}\rightarrow
\underline {\m\mbox{-}\mathcal Gproj}$ given by $Y\mapsto
\left(\begin{smallmatrix}
P\\
Y
\end{smallmatrix}\right)_{\sigma}$, where $P$ is a projective $A$-module such that there is an exact sequence
$0\rightarrow M\otimes Y \stackrel \sigma \rightarrow P \rightarrow
{\rm \Cok} \sigma\rightarrow 0$ with ${\rm \Cok} \sigma\in
A\mbox{-}\mathcal Gproj$.
\end{lem}

\noindent{\bf Proof.} \ Let $_BY$ be Gorenstein-projective. Then
$M\otimes Y$ is Gorenstein-projective, and hence there is an exact
sequence $0\rightarrow M\otimes Y \stackrel \sigma \rightarrow P
\rightarrow {\rm \Cok} \sigma\rightarrow 0$ with $P$ projective and
${\rm \Cok} \sigma\in A\mbox{-}\mathcal Gproj$. Let $g: Y\rightarrow
Y'$ be a $B$-map with $Y, \ Y'\in B\mbox{-}\mathcal Gproj$, and $P'$
a projective $A$-module such that $0\rightarrow M\otimes Y'
\stackrel {\sigma'} \rightarrow P' \rightarrow {\rm \Cok}
\sigma'\rightarrow 0$ is exact with ${\rm \Cok} \sigma'\in
A\mbox{-}\mathcal Gproj$. Since projective $A$-modules are injective
objects in $A\mbox{-}\mathcal Gproj$, it follows that there is a
commutative diagram
\[\xymatrix {
0\ar[r] & M\otimes
Y\ar[d]_{1\otimes g}\ar[r]^-{\sigma} & P\ar[d]_{f}\ar[r]^-{\pi} & {\rm \Cok}\sigma \ar[d]\ar[r] & 0\\
0\ar[r] & M\otimes Y'\ar[r]^-{\sigma'} & P'\ar[r] & {\rm
\Cok}\sigma' \ar[r] & 0.}\] Taking $g={\rm Id}$ we see
$\left(\begin{smallmatrix}
P\\
Y
\end{smallmatrix}\right)_\sigma \cong \left(\begin{smallmatrix}
P'\\
Y
\end{smallmatrix}\right)_{\sigma'}$ in $\underline {\m\mbox{-}\mathcal
Gproj}$. If we have another map $f': P\rightarrow P'$ such that
$f'\sigma = \sigma'(1\otimes g)$, then $f-f'$ factors through
$\Cok\sigma$. Since ${\rm \Cok} \sigma\in A\mbox{-}\mathcal Gproj$,
we have a monomorphism ${\tilde{\sigma}}: {\rm \Cok} \sigma
\rightarrow \tilde{P}$ with $\tilde{P}$ projective. Then we easily
see that $\left(\begin{smallmatrix}
f\\
g
\end{smallmatrix}\right)-\left(\begin{smallmatrix}
f'\\
g
\end{smallmatrix}\right)$ factors through projective $\m$-module $\left(\begin{smallmatrix}
\tilde{P}\\
0
\end{smallmatrix}\right)$, and hence $\underline{\left(\begin{smallmatrix}
f\\
g
\end{smallmatrix}\right)} = \underline{\left(\begin{smallmatrix}
f'\\
g
\end{smallmatrix}\right)}.$ Thus we get a unique functor $j_*: \ {B\mbox{-}\mathcal
Gproj}\rightarrow \underline {\m\mbox{-}\mathcal Gproj}$ given by
$Y\mapsto \left(\begin{smallmatrix}
P\\
Y
\end{smallmatrix}\right)_{\sigma}$ and $g\mapsto \underline{\left(\begin{smallmatrix}
f\\
g
\end{smallmatrix}\right)}$.

\vskip10pt

Assume that $g: Y\rightarrow Y'$ factors through a projective module
$_BQ$ with $g = g_2g_1$. Since $M\otimes Q$ is projective and hence
an injective object in $\underline {A\mbox{-}\mathcal Gproj}$, there
is an $A$-map $\alpha: P\rightarrow M\otimes Q$ such that $1\otimes
g_1 = \alpha \sigma.$ Since $(f-\sigma'(1\otimes g_2)\alpha)\sigma =
0$, there is an $A$-map $\tilde{f}: {\rm \Cok} \sigma \rightarrow
P'$ such that $\tilde{f}\pi = f-\sigma'(1\otimes g_2)\alpha.$ Let
${\tilde{\sigma}}: {\rm \Cok} \sigma \rightarrow \tilde{P}$ be a
monomorphism with $\tilde{P}$ projective. Then we get an $A$-map
$\beta: \tilde{P}\rightarrow P'$ such that $\tilde{f} = \beta
\tilde{\sigma}.$ Thus $\underline{\left(\begin{smallmatrix}
f\\
g
\end{smallmatrix}\right)}$ factors through projective $\m$-module
$\left(\begin{smallmatrix}
M\otimes Q\\
Q
\end{smallmatrix}\right)\oplus \left(\begin{smallmatrix}
\tilde{P}\\
0
\end{smallmatrix}\right)$ with $\underline{\left(\begin{smallmatrix}
f\\
g
\end{smallmatrix}\right)} = \underline{\left(\begin{smallmatrix}
(\sigma'(1\otimes g_2),  \beta)\\
g_2
\end{smallmatrix}\right)} \ \underline{\left(\begin{smallmatrix}
\left(\begin{smallmatrix}
\alpha\\
\tilde{\sigma}\pi
\end{smallmatrix}\right)\\
g_1
\end{smallmatrix}\right)}$. Therefore $j_*: {B\mbox{-}\mathcal
Gproj}\rightarrow \underline {\m\mbox{-}\mathcal Gproj}$ induces a
functor ${ B\mbox{-}\mathcal Gproj}\rightarrow \underline
{\m\mbox{-}\mathcal Gproj},$ again denoted by $j_*$, which is given
by $Y\mapsto \left(\begin{smallmatrix}
P\\
Y
\end{smallmatrix}\right)_{\sigma}$ and $\underline g\mapsto \underline{\left(\begin{smallmatrix}
f\\
g
\end{smallmatrix}\right)}$.

\vskip10pt

By the above argument we know that $j_*$ is full. If
${\left(\begin{smallmatrix}
f\\
g
\end{smallmatrix}\right)}$ factors through projective $\m$-module
$\left(\begin{smallmatrix}
M\otimes Q\\
Q
\end{smallmatrix}\right)\oplus \left(\begin{smallmatrix}
\tilde{P}\\
0
\end{smallmatrix}\right)$, then $g$ factors through
projective module $_BQ$. Thus $j_*$ is faithful. $\s$

\vskip10pt

\subsection{} \ Let $\mathcal{A}$ be a Frobenius category and $\underline{\mathcal{A}}$
the corresponding stable category. Then $\underline{\mathcal{A}}$ is
a triangulated category with shift functor $[1]$ given by $X[1] =
{\Cok} (X\longrightarrow I(X))$, where $I(X)$ is a
projective-injective object of $\mathcal{A}$;  each exact sequence
$0 \rightarrow X \stackrel{u}{\rightarrow} Y
\stackrel{v}{\rightarrow} Z \rightarrow 0$ in $\mathcal{A}$ gives
rise to a distinguished triangle $X
\stackrel{\underline{u}}{\longrightarrow} Y
\stackrel{\underline{v}}{\longrightarrow} Z {\rightarrow}$ \ in
$\underline{\mathcal{A}}$, and each distinguished triangle in
$\underline{\mathcal{A}}$ is of this form up to an isomorphism. See
D. Happel [H], Chapter 1, Section 2. It follows that we have

\vskip10pt

\begin{lem} \ All the functors $i^*, \ i_*, \ i_!, \ j_!, \ j^*, \ j_*$
constructed above are exact functors; and $i_*,  \ j_!$, and $j_*$
are fully faithful.
\end{lem}

\subsection{\bf Proof of Theorem 3.3.} \  By construction ${\rm Ker}j^* =
\{\left(\begin{smallmatrix}
X\\
Q
\end{smallmatrix}\right)_\phi\in \underline {\m\mbox{-}\mathcal
Gproj} \ | \ _BQ \ \mbox{is projective} \}$. By Proposition 3.2
there is an exact sequence $0\rightarrow M\otimes Q\stackrel
\phi\rightarrow X\rightarrow {\rm \Cok}\phi\rightarrow 0$ in
$\m\mbox{-}\mathcal Gproj$. Since $M\otimes Q$ is a projective
$A$-module, and hence an injective object in $\m\mbox{-}\mathcal
Gproj$, it follows that $\phi$ splits and then
$\left(\begin{smallmatrix}
X\\
Q
\end{smallmatrix}\right)_\phi \cong \left(\begin{smallmatrix}
M\otimes Q\\
Q
\end{smallmatrix}\right)_{\rm Id} \oplus \left(\begin{smallmatrix}
X'\\
0
\end{smallmatrix}\right) = \left(\begin{smallmatrix}
X'\\
0
\end{smallmatrix}\right)$ in $\underline {\m\mbox{-}\mathcal
Gproj}$. Thus ${\rm Im}i_* ={\rm Ker}j^*.$

\vskip10pt

In the following $\left(\begin{smallmatrix}
X\\
Y
\end{smallmatrix}\right)_\phi\in \underline {\m\mbox{-}\mathcal
Gproj}$, \ $X'\in \underline {A\mbox{-}\mathcal Gproj}$, and $Y'\in
\underline {B\mbox{-}\mathcal Gproj}$.

\vskip5pt

It is easy to see that a $\m$-map $\left(\begin{smallmatrix}
f\\
0
\end{smallmatrix}\right): \left(\begin{smallmatrix}
X\\
Y
\end{smallmatrix}\right)_\phi\rightarrow \left(\begin{smallmatrix}
X'\\
0
\end{smallmatrix}\right)$ factors through a projective
$\m$-module if and only if the induced $A$-map ${\rm
\Cok}\phi\rightarrow X'$ factors through a projective $A$-module.
This implies that the isomorphism $(2.7)$ induces the following
isomorphism, which are natural in both positions
$$\Hom_{\underline {\m\mbox{-}\mathcal
Gproj}}(\left(\begin{smallmatrix}
X\\
Y
\end{smallmatrix}\right)_\phi, \left(\begin{smallmatrix}
X'\\
0
\end{smallmatrix}\right)) \cong \Hom_{\underline {A\mbox{-}\mathcal
Gproj}} (\Cok\phi, X'), $$ i.e., $(i^\ast, i_\ast)$ is an adjoint
pair.

\vskip5pt

It is easy to see that a $\m$-map $\left(\begin{smallmatrix}
f\\
0
\end{smallmatrix}\right): \left(\begin{smallmatrix}
X'\\
0
\end{smallmatrix}\right)\rightarrow \left(\begin{smallmatrix}
X\\
Y
\end{smallmatrix}\right)_{\phi}$ factors through a projective
$\m$-module if and only if $f: X'\rightarrow X$ factors through a
projective $A$-module. This implies that the isomorphism $(2.8)$
induces the following isomorphism, which are natural in both
positions
$$\Hom_{\underline {\m\mbox{-}\mathcal
Gproj}}(\left(\begin{smallmatrix}
X'\\
0
\end{smallmatrix}\right), \left(\begin{smallmatrix}
X\\
Y
\end{smallmatrix}\right)_\phi)
\cong \Hom_{\underline {A\mbox{-}\mathcal Gproj}} (X', X),
$$ i.e., $(i_\ast, i^!)$ is an adjoint pair.

\vskip5pt

Note that $M\otimes Q$ is a projective $A$-module for any projective
$B$-module $Q$. It is easy to see that a $\m$-map
$\left(\begin{smallmatrix}
\phi({\rm Id_M}\otimes g)\\
g
\end{smallmatrix}\right): \left(\begin{smallmatrix}
M\otimes Y'\\
Y'
\end{smallmatrix}\right)_{\rm Id}\rightarrow \left(\begin{smallmatrix}
X\\
Y
\end{smallmatrix}\right)_{\phi}$ factors through a projective
$\m$-module if and only if $g: Y'\rightarrow Y$ factors through a
projective $B$-module. This implies that the isomorphism $(2.9)$
induces the following isomorphism, which are natural in both
positions
$$\Hom_{\underline {\m\mbox{-}\mathcal
Gproj}}(\left(\begin{smallmatrix}
M\otimes Y'\\
Y'
\end{smallmatrix}\right)_{\rm Id}, \left(\begin{smallmatrix}
X\\
Y
\end{smallmatrix}\right)_{\phi}) \cong \Hom_{\underline {B\mbox{-}\mathcal
Gproj}} (Y', Y), $$ i.e., $(j_!, j^\ast)$ is an adjoint pair.

\vskip5pt

Let $\left(\begin{smallmatrix}
f\\
g
\end{smallmatrix}\right): \left(\begin{smallmatrix}
X\\
Y
\end{smallmatrix}\right)_\phi\rightarrow \left(\begin{smallmatrix}
P'\\
Y'
\end{smallmatrix}\right)_{\sigma}$ be a $\m$-map,  $0\rightarrow M\otimes Y' \stackrel \sigma \rightarrow P'
\rightarrow {\rm \Cok} \sigma\rightarrow 0$ an exact sequence with
$P'$ projective and ${\rm \Cok} \sigma\in A\mbox{-}\mathcal Gproj$.
In the proof of Lemma 3.4 we know that $\left(\begin{smallmatrix}
f\\
g
\end{smallmatrix}\right)$ factors
through a projective $\m$-module if and only if $g: Y\rightarrow Y'$
factors through a projective $B$-module. This implies that the map
$\underline g\mapsto  \underline{\left(\begin{smallmatrix}
f\\
g
\end{smallmatrix}\right)}$ gives rise to the following isomorphism, which is natural in both positions
$$\Hom_{\underline
{\m\mbox{-}\mathcal Gproj}}(\left(\begin{smallmatrix}
X\\
Y
\end{smallmatrix}\right)_{\phi}, \left(\begin{smallmatrix}
P'\\
Y'
\end{smallmatrix}\right)) \cong \Hom_{\underline {\m\mbox{-}\mathcal
Gproj}} (Y, Y'),$$ i.e., $(j^\ast, j_\ast)$ is an adjoint pair. Now
the first assertion follows from Lemmas 3.5 and 1.3.

\vskip10pt

Assume that $A$ and $B$  are in additional finite-dimensional
algebras over a field $k$. Note that $\m\mbox{-}\mathcal Gproj$ is a
resolving subcategory of $\m$-mod (see e.g. Theorem 2.5 in [Hol]).
Since $\m$ is a Gorenstein algebra, it is well-known that
$\m\mbox{-}\mathcal Gproj$ contravariantly finite in $\m$-mod (see
Theorem 11.5.1 in [EJ], where the result is stated for arbitrary
$\m$-modules, but the proof holds also for finitely generated
modules. See also Theorem 2.10 in [Hol]). Then by Corollary 0.3 of
H. Krause and \O. Solberg [KS], which asserts that a resolving
contravariantly finite subcategory in $A$-mod is also covariantly
finite in $A$-mod, $\m\mbox{-}\mathcal Gproj$ is functorially finite
in $A$-mod, and hence $\m\mbox{-}\mathcal Gproj$ has
Auslander-Reiten sequences, by Theorem 2.4 of M. Auslander and S. O.
Smal${\o}$ [AS]. Since each distinguished triangle in the stable
category $\underline{\mathcal A}$ of a Frobebius category $\mathcal
A$ is induced by an exact sequence in $\mathcal A$, $\underline
{\m\mbox{-}\mathcal Gproj}$ has Auslander-Reiten triangles. By
assumption $\m$ is finite-dimensional $k$-algebra, thus $\underline
{\m\mbox{-}\mathcal Gproj}$ is a Hom-finite $k$-linear Krull-Schmidt
category, and hence by Theorem I.2.4 of I. Reiten and M. Van den
Bergh [RV] $\underline {\m\mbox{-}\mathcal Gproj}$ has a Serre
functor. Now the second assertion follows from Theorem 7 of P.
J$\o$rgensen [J], which claims that any recollement of a
triangulated category with a Serre functor is symmetric. $\s$

\vskip10pt

\subsection{} \ By Theorem 3.3 we have

\vskip10pt

\begin{cor} \ Let $A$ be a Gorenstein algebra, and $T_2(A)
= \left(\begin{smallmatrix}
A&A\\
0&A
\end{smallmatrix}\right)$.
Then we have a recollement of triangulated categories

\begin{center}
\begin{picture}(120,45)
\put(-5,20){\makebox(-16,1){$\underline {A\mbox{-}\mathcal Gproj}$}}
\put(40,30){\vector(-1,0){30}} \put(10,20){\vector(1,0){30}}
\put(40,10){\vector(-1,0){30}} \put(70,20) {\makebox(16,1)
{$\underline {T_2(A)\mbox{-}\mathcal Gproj}$}}
\put(145,30){\vector(-1,0){30}} \put(115,20){\vector(1,0){30}}
\put(145,10){\vector(-1,0){30}} \put(165,20){\makebox(18, 0.5)
{$\underline {A\mbox{-}\mathcal Gproj}$\ ;}}
\end{picture}
\end{center}
and it is symmetric if $A$ and $B$ are finite-dimensional algebras
over a field.\end{cor}

\vskip10pt

For the first part of Corollary 3.6 see also Theorem 3.8 in [IKM].

\end{document}